\documentclass[11pt, oneside]{article}   	
\usepackage{geometry}                		
\usepackage{graphicx}				
\usepackage{amssymb,amsthm,amsmath}

\title{Solution to a conjecture on edge rings with 2-linear resolutions}
\author{Ralf Fr\"oberg}
\date{}							

\begin{document}
\newtheorem{theorem}{Theorem}

\newtheorem{lemma}[theorem]{Lemma}

\newtheorem{cor}{Corollary}

\newtheorem*{conj}{Conjecture}

\newtheorem{definition}[theorem]{Definition}

\newtheorem{rmk}{Remark}

\maketitle
\begin{abstract}
For a graph $G=(V,E)$ the edge ring $k[G]$ is $k[x_1,\ldots,x_n]/I(G)$, where $n=|V|$ and $I(G)$ is generated by $\{ x_ix_j;\{ i,j\}\in E\}$.
The conjecture we treat is the following. 

\medskip
\begin{conj} If $k[G]$ has a 2-linear resolution, then the projective dimension of $K[G]$, pd$(k[G])$, equals the maximal degree of a vertex in $G$.
\end{conj}

As far as we know,
this conjecture is first mentioned in a paper by Gitler and Valencia, \cite[Conjecture 4.13]{g-v}, and there it is called the Eliahou-Villarreal conjecture. 
The conjecture is treated in a recent paper by Ahmed, Mafi, and Namiq, \cite{a-m-n}. That there are counterexamples was noted already by
Moradi and Kiani, \cite{k-m}. By interpreting $k[G]$ as a Stanley-Reisner ring, we are able
to characterize those graphs for which the conjecture holds.
\end{abstract}

\section{Introduction}
An edge ring $k[G]$ has a 2-linear resolution if and only if the complement graph $\overline G$ is chordal, \cite[Theorem 1]{f1}, i.e., every cycle of length $\ge3$ in 
$\overline G$ has a chord. (The complement graph $\overline G$ has the same vertex set as $G$, and $\{ i,j\}$ is an edge in $\overline G$ if and only
if $\{ i,j\}$ is not an edge in $G$.) There are other proofs of the theorem in \cite[Theorem 2.1]{e-g}, \cite[[Theorem 3.4]{d-e}, \cite[Corollary 3.3]{n}, and \cite[Theorem 9.2.12]{h-h}.

However, we prefer to use the characterization of edge rings with 2-linear resolution in terms of Stanley-Reisner rings.
The Stanley-Reisner ring of a simplicial complex $\Delta$ with vertex set $\{ 1,2,\ldots,n\}$ is $k[x_1,\ldots,x_n]/I(\Delta)$, where $I(\Delta)$ is
generated by those squarefree monomials $x_{i_1}\cdots x_{i_k}$ for which
$\{ i_1,\ldots,i_k\}$ is not a face in $\Delta$. Then $k[\Delta]$ has a 2-linear resolution if and only if $\Delta$ is a quasi-forest, \cite{f1,f2}. 
Quasi-forests are also named fat forests or generalized $d$-trees, but we have chosen the name quasi-forest here, because it seems to be the most common. A quasi-forest
is defined recursively in the following way. A simplex $F_1$ is a quasi-forest. If $F_1\cup F_2\cup\cdots\cup F_{i-1}$, $F_j$ simplices for $1\le j\le i-1$, is a 
qusi-forest, then, if $F_i$ is a simplex, $F_1\cup F_2\cup\cdots\cup F_i$ is a quasi-forest if $G_i=(F_1\cup F_2\cup\cdots\cup F_{i-1})\cap F_i$ is a simplex.
If $F_1\cup F_2\cup\cdots\cup F_{i-1}$ and $F_i$ are disjoint, then $G_i=\emptyset$ of dimension -1.

For a simplicial complex $\Delta$, we let $G(\Delta)$ be its 1-skeleton. For a graph $G$ we let $\Delta(G)$ be the largest simplicial complex with
$G$ as 1-skeleton. If $k[F]=k[x_1,\ldots,x_n]/I(F)$ is the Stanley-Reisner ring of the simplicial complex $F$, then $I(F)$ is generated in degree 2
if and only if $\Delta(G(F))=F$. We can go back and forth between Stanley-Reisner rings with 2-linear resolution and edge rings with 2-linear resolution since 
$k[G]$ has a 2-linear resolution if and only if $k[\Delta(\overline G)]$ has a 2-linear resolution, and $k[\Delta]$ has a 2-linear resolution if and only if
$k[\overline{G(\Delta)}]$ has a 2-linear resolution.
Let $F=F_1\cup\cdots\cup F_k$ be a fat forest, and $H$ the 1-skeleton of $F$. Then the edge ring of $\overline H$ is the Stanley-Reisner ring of $F$,
and thus has a 2-linear resolution. 

Suppose $\dim F_i=d_i-1$ and $\dim G_i=r_i-1$. If $d_1\ge d_2\ge\cdots\ge d_i$ and $d_j=r_j+1$ for all $j$, then the 1-skeleton $H$ 
is what is called a $(d_1,d_2,\ldots,d_i)$-tree in \cite{a-m-n}. 

A vertex $v$ in a simplicial complex $\Delta$ is called free if $v$ belongs to a unique facet (maximal face). 

\section{Hilbert series and Betti numbers of fat forests}
If $k[x_1,\ldots,x_n]/I=S/I$ has a 2-linear resolution it looks like this:

$$0\leftarrow S/I\leftarrow S\leftarrow S[-2]^{b_1}\leftarrow S[-3]^{b_2}\leftarrow\cdots\leftarrow S[-p-1]^{b_p}\leftarrow 0$$
where $S[-i]$ means that we have shifted degrees of $S$ $i$ steps.

If we restrict the sequence to a certain degree, we get an exact sequence of vector spaces, thus the alternating sum of dimensions is 0.
Using this we get that the
Hilbert series $\sum_{i\ge0}\dim_kk[\Delta]_it^i$ of $k[\Delta]$ with 2-linear resolution equals 
$\frac{1-\beta_{1,2}t^2+\beta_{2,3}t^3-\cdots(-1)^{p}\beta_{p,p+1}t^{p+1}}{(1-t)^n}$, where $\beta_{i,j}$ are 
the graded Betti numbers $\dim_k{\rm Tor}^S_{i,j}(k[\Delta],k)$, and $n$ is the
number of vertices in $\Delta$.
If one is interested in the Betti numbers $\beta_{i,j}=\dim_k{\rm Tor}^S_{i,j}(S/I,k)$ of
Stanley-Reisner ring of a fat tree, since it has 2-linear resolutionsons, one could just as well study its Hilbert series, because this
contains the same information as the set of Betti numbers. We now give the main result from \cite{f2}. To make this paper self contained
we repeat the short proof.

\begin{theorem}\label{y} Let $F=F_1\cup\cdots\cup F_k$ be a quasi-tree with $F_i$ a simplex of dimension $d_i$
 and $(F_1\cup\cdots\cup F_{j-1})\cap F_j$ a simplex of dimension $r_j$. 
 Then the Hilbert series of $k[F]$
 is $\sum_{i=1}^k\frac{1}{(1-t)^{d_i+1}}-\sum_{i=2}^{k}\frac{1}{(1-t)^{r_i+1}}$.
 The projective dimension is $\sum_{i=1}^kd_i-\sum_{i=2}^kr_i+1-\min\{ r_i\}-2$. 
 The depth of $k[F]$ is $\min\{ r_i\}+2$, and $F$ is CM (Cohen-Macaulay) 
 if and only if there is a $d$ such that
 $d_i=d$ for all $i$ and $r_i=d-1$ for all $i$.
  \end{theorem}
 
 \noindent
 {\em Proof} A $d$-simplex has Hilbert series $\frac{1}{(1-t)^{d+1}}$. We have to subtract the Hilbert series of the $G_i$'s 
 from $\sum\frac{1}{(1-t)^{d_i+1}}$, because they are counted twice in $\sum_{i=1}^k\frac{1}{(1-t)^{d_i+1}}$.
The number of vertices of $F$ is $\sum_{i=1}^k(d_i+1)-\sum_{i=2}^k(r_i+1)=
 \sum_{i=1}^kd_i-\sum_{i=2}^kr_i+1=n$, so the degree of the numerator $p(t)$ of the Hilbert series
 $\frac{p(t)}{(1-t)^n}$ of $k[F]$ is $n-\min\{ r_i\}-1$ so the projective dimension is $n-\min\{ r_i\}-2$,
 and the depth of $k[F]$ is $\min\{ r_i\}+2$ by the Auslander-Buchsbaum formula. We have
 $\dim k[F]=1+\max\{ d_i\}$, depth $k[F]=\min\{ r_i\}+2$, and $d_i>r_i$ for all $i$. The only possibility
 for $\dim k[F]={\rm depth}\, k[F]$ is that there is a $d$ such that $d_i=r_i+1=d$ for all $i$.
 
 \section{Main result}
 We now come to the main result, the characterization of edge rings with 2-linear resolution for which the conjecture holds.
 We first state it in terms of Stanley-Reisner rings $k[\Delta]$.
 If $\Delta$ is a simplex, the conjecture trivially holds. Thus we may suppose that $\Delta$ has at least two facets.
 
 \begin{theorem}\label{z}
 Suppose $\Delta$ is a quasi-forest with at least two facets, in particular $k[\Delta]$ has a 2-linear resolution. 
 Let $r=\min\{ r_i\}$. Then the conjecture holds if and only if there is a free vertex $v$ in a facet $F$ of $\Delta$, $\dim F=r+1$, 
 $F$ connected to the remaning part of $\Delta$ in a simplex of dimension $r$.
 \end{theorem}
 
 \noindent
 {\em Proof} Let $v$ be a vertex in $\Delta$. Suppose $v$ belongs to $F_{i_1},\ldots,F_{i_k}$. The dimension of $F_{i_1}$ is at least
 $r+1$ so $F_{i_1}$ has at least $r+2$ vertices. There must be at least one vertex in each $F_{i_j}$, $j\ge2$,
 which is not in the union of the others. Thus there are $n-(r+2+k-1)$ vertices which are not neighbours. (That $u$ is not a neighbour
 to $v$ means exactly that $u$ is a neighbour to $v$ in $\overline{G(\Delta)}$.) The projective dimension is $n-r-2$. 
 If the numbers are equal, we get $k=1$, and furthermore that $F_{i_1}$ has precisely $r+2$ vertices, so $d_{i_1}=r+1$.

\begin{cor}\label{x}
Let $G$ be a graph and suppose that $k[G]$ has a 2-linear resolution. Then the conjecture holds for $G$ if and only
$\overline G$ is the 1-skeleton of a simplicial complex as in Theorem \ref{z}.
\end{cor}

\begin{cor}
The smallest counterexample to the conjecture is the the edge ring of a 4-cycle, noted already in \cite{k-m}.
This corresponds to the Stanley-Reisner ring on four vertices $\{1,2,3,4\}$ with facets $\{1,2\},\{3,4\}$. 
\end{cor}

\noindent
{\em Proof} If $G=C_4$, then $\Delta(\overline G)$ is the complex $F$ on four vertices $\{1,2,3,4\}$ with facets $\{1,2\},\{3,4\}$. 
We have pd$(k[C_4])=$pd$(k[F])=4-(-1)-2=3$ and $\max\deg\{v;v\in G\}=2$.

\begin{cor}
If the quasi-forest $\Delta$ has an isolated vertex $v$, then the conjecture holds.
\end{cor}

\noindent
{\em Proof} We have $\dim\{v\}=0$ and $r=-1$.

\begin{rmk}
This is equivalent to \cite[Theorem 4.11]{a-m-n}, which considers the case when there is a vertex in the graph with all other vertices a neighbours (a full vertex).
\end{rmk}

\begin{theorem}\label{u}\cite[Theorem 4.21]{a-m-n}
Let $G$ be a graph such that $\overline G$ is a $(d_1,d_2,\ldots,d_q)$-tree. Then the conjecture holds.
\end{theorem}

\noindent
{\em Proof}
This follows from Corollary \ref{x}.

\begin{cor}
If $k[G]$ is a Cohen-Macaulay ring with 2-linear resolution, then the conjecture holds for $k[G]$.
\end{cor}

\noindent
{\em Proof} If $k[G]$ is Cohen-Macaulay, then $k[G]=k[\Delta]$, and the complement of the1-skeleton of $\Delta$ is a $(d,d,\ldots,d)$-tree for some $d$ by 
definition of $(d_1,\ldots,d_q)$-trees and by Theorem \ref{y}. 
Thus the result follows Theorem \ref{u}.

\begin{rmk}
The authors of \cite{a-m-n} also give an example to show that the difference pd$(k[G])-\max\deg \{ v;v\in G\}$ can be arbirarily large. The example is the 
the complement of the $r$-barbell, i.e., the complement to the graph consisting of two complete graphs $K_r$ with a bridge. We give a slightly smaller example.
\end{rmk}

\begin{cor}
Let $K_{r,r}$ be the complete bipartite graph. Then pd$(k[K_{r,r}])-\max\deg\{v;v\in K_{r,r}\}=r$.
\end{cor}

\noindent
{\em Proof} $K_{r,r}$ is the complement to the 1-skeleton of the disjoint union $S_{r-1}\sqcup S_{r-1}$ of two $(r-1)$-simplices $S_{r-1}$.
Now the edge ring $k[K_{r,r}]$ equals the Stanley-reisner ring $k[S_{r-1}\sqcup S_{r-1}]$ and pd$(k[S_{r-1}\sqcup S_{r-1}])-\max\deg\{v;v\in K_{r,r}\}=2r-(-1)-2-r=r-1$.

\begin{rmk}
In \cite{a-m-n} it is claimed that the conjectrure is stated in \cite{e-v}. This seems to be a mistake, the conjecture is not mentioned in that paper, and the name of the
conjecture is a bit of the mystery.
\end{rmk}

\end{document}